\documentclass[12pt]{amsart}
\usepackage{amscd}
\def\frk{\frak}               

\def\Q! ! q{{\frk Q}}

\def\Phi{{\frk n}}
\def\Phi{{\frk N}}
\def\opn#1#2{\def#1{\operatorname{#2}}} 
\opn\chara{char}
\opn\length{\ell}
\opn\pd{pd}
\opn\rk{rk}
\opn\projdim{proj\,dim}
\opn\injdim{inj\,dim}
\opn\rank{rank}
\opn\depth{depth}
\opn\grade{grade}
\opn\height{height}
\opn\embdim{emb\,dim}
\opn\codim{codim}

\opn\Tr{Tr}
\opn\bigrank{big\,rank}
\opn\superheight{superheight}\opn\lcm{lcm}
\opn\trdeg{tr\,deg}%
\opn\reg{reg}
\opn\lreg{lreg}
\opn\ini{in}
\opn\lpd{lpd}
\opn\div{div}
\opn\Div{Div}
\opn\cl{cl}
\opn\Cl{Cl}
\opn\Spec{Spec}
\opn\Supp{Supp}
\opn\supp{supp}
\opn\Sing{Sing}
\opn\Ass{Ass}
\opn\Ann{Ann}
\opn\Rad{Rad}
\opn\Soc{Soc}
\opn\Im{Im}
\opn\Ker{Ker}
\opn\Coker{Coker}
\opn\Am{Am}
\opn\Hom{Hom}
\opn\Tor{Tor}
\opn\Ext{Ext}
\opn\End{End}
\opn\Aut{Aut}
\opn\id{id}

\opn\nat{nat}
\opn\pff{pf}
\opn\Pf{Pf}
\opn\GL{GL}
\opn\SL{SL}
\opn\mod{mod}
\opn\ord{ord}
\opn\Gin{Gin}
\opn\Hilb{Hilb}
\opn\aff{aff}
\opn\con{conv}
\opn\relint{relint}
\opn\st{st}
\opn\lk{lk}
\opn\cn{cn}
\opn\core{core}
\opn\vol{vol}
\opn\dist{dist}
\opn\link{link}
\opn\star{star}
\opn\gr{gr}

\def\pot#1#2{#1[\kern-0.28ex[#2]\kern-0.28ex]}

\opn\dirlim{\underrightarrow{\lim}}
\opn\inivlim{\underleftarrow{\lim}}
\let\Union=\bigcup

\def\Implies{\ifmmode\Longrightarrow \else
        \unskip${}\Longrightarrow{}$\ignorespaces\fi}
\def\implies{\ifmmode\Rightarrow \else
        \unskip${}\Rightarrow{}$\ignorespaces\fi}
\def\iff{\ifmmode\Longleftrightarrow \else
        \unskip${}\Longleftrightarrow{}$\ignorespaces\fi}

\let\:=\colon
\newtheorem{Theorem}{Theorem}[section]
\newtheorem{Lemma}[Theorem]{Lemma}

\newtheorem{Example}[Theorem]{Example}

\let\epsilon\varepsilon
\let\phi=\varphi
\let\kappa=\varkappa
\def\qed{\ifhmode\textqed\fi
      \ifmmode\ifinner\quad\qedsymbol\else\dispqed\fi\fi}
\def\textqed{\unskip\nobreak\penalty50
       \hskip2em\hbox{}\nobreak\hfil\qedsymbol
       \parfillskip=0pt \finalhyphendemerits=0}
\def\dispqed{\rlap{\qquad\qedsymbol}}

\opn\dis{dis}
\def\pnt{{\raise0.5mm\hbox{\large\bf.}}}

\opn\Lex{Lex}
\begin{document}
\title{Cohen--Macaulay polymatroidal ideals}
\author{J\"urgen Herzog and Takayuki Hibi}
\subjclass{13C14, 05B35}
\thanks{This paper was completed while the authors stayed at the Mathematisches Forschungsinstitut in Oberwolfach in the frame of the Research in Pairs Program}
\address{J\"urgen Herzog, Fachbereich Mathematik und
Informatik, Universit\"at Duisburg-Essen, Campus Essen,
45117 Essen, Germany}
\email{juergen.herzog@uni-essen.de}
\address{Takayuki Hibi, Department of Pure and Applied Mathematics,
Graduate School of Information Science and Technology,
Osaka University, Toyonaka, Osaka 560-0043, Japan}
\email{hibi@math.sci.osaka-u.ac.jp}

\maketitle
\begin{abstract}
All Cohen--Macaulay polymatroidal ideals are classified.
The Cohen--Macaulay polymatroidal ideals 
are precisely the 
principal ideals, the 
Veronese ideals, and the 
squarefree Veronese ideals.
\end{abstract}

\section*{Introduction}
 Our goal  is to classify all
Cohen--Macaulay polymatroid ideals.
It can be expected that such classification would be possible.
Because, it seems likely that Cohen--Macaulay monomial ideals 
with linear resolutions are quite rare and it is known that 
every polymatroid ideal has a linear resolution.
Our main result says that
a polymatroidal ideal $I$ is Cohen--Macaulay 
if and only if
$I$ is 
a principal ideal,  a Veronese ideal, or
a squarefree Veronese ideal,
see Theorem \ref{polymatroidal}.

\section{Monomial ideals with linear quotients}
Let $K$ be a field and $S = K[x_1, \ldots, x_n]$
the polynomial ring in $n$ variables over $K$
with each $\deg x_i = 1$.
Let $I \subset S$ be a monomial ideal and
$G(I)$ its unique minimal monomial generators.

A {\em vertex cover} of $I$ is a subset $W$ of
$\{ x_1, \ldots, x_n \}$ such that each $u \in G(I)$ 
is divided by some $x_i \in W$.  
Such a vertex cover $W$ is called {\em minimal}
if no proper subset of $W$ is a vertex cover of $I$.

A monomial ideal is called {\em unmixed}
if all minimal vertex covers of $I$ have the same cardinality.
If $I$ is {\em Cohen--Macaulay}, i.e., the quotient ring $S/I$
is Cohen--Macaulay, then $I$ is unmixed.
Let $h(I)$ denote the minimal cardinality of the vertex covers
of $I$.  It then follows that
\begin{eqnarray}
\label{dim}
\dim S/I = n - h(I).
\end{eqnarray}

We say that a monomial ideal $I \subset S$ has {\em linear 
quotients} if there is an ordering $u_1, \ldots, u_s$ of 
the monomials belonging to $G(I)$ with
$\deg u_1 \leq \deg u_2 \leq \cdots\leq \deg u_s$
such that,
for each $2 \leq j \leq s$,
the colon ideal 
$(u_{1}, u_{2}, \ldots, u_{j-1}) : u_{j}$
is generated by a subset of $\{ x_1, \ldots, x_ n \}$.

It is known, e.g., \cite[Lemma 4.1]{ConcaHerzog} that 
if a monomial ideal $I$ generated in one degree
has linear quotients, then $I$ has a linear resolution.

Let $I$ be a monomial ideal with linear quotient
with respect to the ordering 
$u_1, \ldots, u_s$ of 
the monomials belonging to $G(I)$.
We write $q_j(I)$ for the number of variables which is 
required to generate the colon ideal 
$(u_{1}, u_{2}, \ldots, u_{j-1}) : u_{j}$.
Let $q(I) = \max_{2 \leq j \leq s} q_j(I)$.
It is proved
\cite[Corollary 1.6]{HerzogTakayama}
that the length of the minimal free resolution of
$S/I$ over $S$ is equal to $q(I) + 1$.
Hence 
\begin{eqnarray}
\label{depth}
\depth S/I = n - q(I) - 1.
\end{eqnarray}
Thus in particular the integer $q(I)$ is independent of
the particular choice of the ordering 
of the monomials which gives  linear quotients.

By using the formulae (\ref{dim}) and (\ref{depth}),
it follows that
a monomial ideal $I$ with linear quotients
is Cohen--Macaulay if and only if
$h(I) = q(I) + 1$.

\section{Review on polymatroidal ideals}
One of the important classes of monomial ideals 
with linear quotients
is the class of polymatroid ideals.

Let, as before, $K$ be a field
and $S = K[x_1, \ldots, x_n]$
the polynomial ring in $n$ variables over $K$
with each $\deg x_i = 1$.
Let $I \subset S$ be a monomial ideal generated in one degree.
We say that $I$ is {\em polymatroidal}
if the following ``exchange condition'' is satisfied:
For monomials $u = x_1^{a_1} \cdots x_n^{a_n}$ and
$v = x_1^{b_1} \cdots x_n^{b_n}$ belonging to $G(I)$
and for each $i$ with $a_i > b_i$, one has $j$ with
$a_j < b_j$ such that $x_j u / x_i \in G(I)$.
The reason why we call such an ideal
polymatroidal is that the monomials of the ideal correspond
to the bases of a discrete polymatroid \cite{HerzogHibi}.
The polymatroidal ideal $I$ is called {\em matroidal} if
$I$ is generated by squarefree monomials.

The exchange property for polymatroidal 
ideals has a ``dual version'' stated below.

\begin{Lemma}
\label{juergen}
Let $I$ be a polymatroidal ideal.  Then,
for monomials $u = x_1^{a_1} \cdots x_n^{a_n}$ and
$v = x_1^{b_1} \cdots x_n^{b_n}$ belonging to $G(I)$
and for each $i$ with $a_i < b_i$, one has $j$ with
$a_j > b_j$ such that $x_i u / x_j \in G(I)$.
\end{Lemma}

\begin{proof}
We introduce the distance of $u$ and $v$ by setting
$\dist(u, v) = \frac{1}{2}\sum_{q=1}^{n} |a_q - b_q|$.
Fix $i$ with $a_i < b_i$.
If there is $k_1 \neq i$ with $a_{k_1} < b_{k_1}$, then
there is $\ell_1$ with $a_{\ell_1} > b_{\ell_1}$
such that $w_1 = x_{\ell_1} v / x_{k_1} \in G(I)$.
Let $w_1 = x_1^{c_1} \cdots x_n^{c_n}$.  Then
$c_i = b_i$ and $\dist(u, w_1) < \dist(u, v)$.
Again, if there is $k_2 \neq i$ with $a_{k_2} < c_{k_2}$,
then there is $\ell_2$ with $a_{\ell_2} > c_{\ell_2}$
such that $w_2 = x_{\ell_2} w_1 / x_{k_2} \in G(I)$.
Let $w_2 = x_1^{d_1} \cdots x_n^{d_n}$.
Then $d_i = b_i$ and $\dist(u, w_2) < \dist(u, w_1)$.
Repeating these procedures yields
$w^* = x_1^{q_1} \cdots x_n^{q_n} \in G(I)$ with
$q_i = b_i > a_i$ and $q_j \leq a_j$ for all $j \neq i$.
One has $j_0 \neq i$ with $q_{j_0} < a_{j_0}$.
Then $x_i u / x_{j_0} \in G(I)$, as desired.\hspace{9cm}
\end{proof}

It is known \cite[Theorem 5.2]{ConcaHerzog} that
a polymatroidal ideal has linear quotients with respect to
the reverse lexicographic order $<_{rev}$ 
induced by the ordering $x_1 > x_2 > \cdots > x_n$. 
More precisely,
if $I$ is a polymatroidal ideal and
if $u_1, \ldots, u_s$ are the monomials belonging to $G(I)$
ordered by the reverse lexicographic order, i.e.,
$u_s <_{rev} \cdots <_{rev} u_2 <_{rev} u_1$, then
the colon ideal $(u_1, \ldots, u_{j-1}) : u_j$
is generated by a subset of $\{ x_1, \ldots, x_n \}$.

The product of polymatroidal ideals is again
polymatroidal (\cite{ConcaHerzog} and \cite{HerzogHibi}).
In particular each power of a polymatroidal ideal is
polymatroidal.

We close the present section with polymatroidal ideals 
of special kinds which are of great interest to us.

\begin{Example}
\label{example}
{\em
(a)
The {\em Veronese ideal} of degree $d$ in the variables
$x_{i_1}, \ldots, x_{i_t}$ is the ideal of $S$
which is generated by all monomials
in $x_{i_1}, \ldots, x_{i_t}$ of degree $d$.
The Veronese ideal is polymatroidal 
and is Cohen--Macaulay.

(b)
The {\em squarefree Veronese ideal} of degree $d$ 
in the variables
$x_{i_1}, \ldots, x_{i_t}$ is the ideal of $S$
which is generated by all squarefree monomials
in $x_{i_1}, \ldots, x_{i_t}$ of degree $d$.
The squarefree Veronese ideal is matroidal 
and is Cohen--Macaulay.
}
\end{Example}

\section{Classification of Cohen-Macaulay polymatroidal ideals}
We now classify all Cohen--Macaulay polymatroidal ideals.
Recall that the support of a monomial
$u = x_1^{a_1} \cdots x_n^{a_n}$
is $\supp(u) = \{ x_i : a_i \neq 0 \}$.

\begin{Lemma}
\label{radical}
If $I \subset S$ is a Cohen--Macaulay polymatroidal ideal, 
then its radical $\sqrt{I}$ is squarefree Veronese.
\end{Lemma}

\begin{proof}
Let $I \subset S$ be a Cohen--Macaulay polymatroidal ideal. We may  assume that  $\Union_{u\in G(I)}\supp(u)=\{x_1,\ldots,x_n\}$.
Let $u \in G(I)$ be a monomial
for which $|\supp(u)|$ is minimal.
Let, say, $\supp(u) = \{ x_{n-d+1}, x_{n-d+2}, \ldots, x_n \}$.
Let $J$ denote the monomial ideal generated by those monomials
$w \in G(I)$ such that $w$ is bigger than $u$ 
with respect to the reverse lexicographic order.
We know that the colon ideal $J : u$
is generated by a subset $M$ of $\{ x_{1}, \ldots, x_{n} \}$.
We claim that $\{ x_{1}, \ldots, x_{n-d} \} \subset M$.
For each $1 \leq i \leq n - d$, 
there is a monomial belonging to $G(I)$
which is divided by $x_i$.
It follows from Lemma \ref{juergen} that
there is a variable $x_j$ with $n-d+1 \leq j \leq n$
such that $v = x_i u / x_j \in G(I)$.
One has $v \in J$.  Since $x_i u = x_j v \in J$, one has
$x_i \in J : u$, as required.
Consequently, one has $q(I) \geq n - d$.
Since $I$ is Cohen--Macaulay, it follows that
$h(I) \geq n - d + 1$.
It then turns out that,
for each subset $W \subset \{ x_1, \ldots, x_n \}$
with $|W| = d$, the set
$\{ x_1, \ldots, x_n \} \setminus W$
cannot be a vertex cover of $I$.
Hence for each subset $W \subset \{ x_1, \ldots, x_n \}$
with $|W| = d$ there is a monomial $w \in G(I)$
with $\supp(w) \subset W$.  Since
$|\supp(w)| \geq |\supp(u)| = d$, one has $\supp(w) = W$.
Hence $\sqrt{I}$ is generated by all squarefree monomials
of degree $d$ in $x_1, \ldots, x_n$.\hspace{2cm}
\end{proof}

\begin{Theorem}
\label{polymatroidal}
A polymatroidal ideal $I$ is Cohen--Macaulay if and only if
$I$ is 
\begin{enumerate}
\item[(i)] a principal ideal, 
\item[(ii)] a Veronese ideal, or
\item[(iii)] a squarefree Veronese ideal.
\end{enumerate}
\end{Theorem}

\begin{proof}
By using Lemma \ref{radical} we assume that $\sqrt{I}$
is generated by all squarefree monomials
of degree $d$ in $x_1, \ldots, x_n$, where
$2 \leq d < n$.
One has $h(I) = h(\sqrt{I}) = n - d + 1$.
Suppose that $I$ is not squarefree
(or, equivalently,
each monomial belonging to $G(I)$ is of degree $> d$).
Let $u = \prod_{i=n-d+1}^{n} x_i^{a_{i}}
\in G(I)$ be a monomial with
$\supp(u) = \{ x_{n-d+1}, x_{n-d+2}, \ldots, x_n \}$.
For a while, we assume that
$(*)$ there is a monomial
$v = \prod_{i=1}^{n} x_i^{b_{i}}\in G(I)$ with
$b_{n-d+1} > a_{n-d+1}$.
Let $J$ denote the monomial ideal generated by those monomials
$w \in G(I)$ such that $w$ is bigger than $u$ with respect to the
reverse
lexicographic order.
As was shown in the proof of Lemma \ref{radical},
the colon ideal $J : u$
is generated by a subset $M$ of $\{ x_{1}, \ldots, x_{n} \}$
with $\{ x_{1}, \ldots, x_{n-d} \} \subset M$.
We claim that $x_{n-d+1} \in J : u$.
By using Lemma \ref{juergen} our assumption $(*)$ guarantees that
there is a variable $x_j$ with $n-d+1 < j \leq n$ such that
$u_0 = x_{n-d+1} u / x_j \in G(I)$.  Since $u_0 \in J$,
one has $x_{n-d+1} \in M$.
Hence $q(I) \geq n - d + 1$.
Thus $h(I) < q(I) + 1$ and $I$ cannot be Cohen--Macaulay.

To complete our proof, we must examine our assumption $(*)$.
For each $d$-element subset
$\sigma = \{ x_{i_1}, x_{i_2}, \ldots, x_{i_d} \}$
of $\{ x_{1}, \ldots, x_{n} \}$, there is a monomial
$u_\sigma \in G(I)$ with
$\supp(u_\sigma) = \sigma$.
If there are $d$-element subset $\sigma$ and $\tau$
of $\{ x_{1}, \ldots, x_{n} \}$ and a variable
$x_{i_0} \in \sigma \cap \tau$ with
$a_{i_0} < b_{i_0}$,
where
$a_{i_0}$ (resp. $b_{i_0}$) is the power of $x_{i_0}$
in $u_\sigma$ (resp. $u_\tau$),
then after relabelling the variables if necessarily
we may assume that
$\sigma = \{ x_{n-d+1}, x_{n-d+2}, \ldots, x_n \}$
with $i_0 = n - d + 1$.  In other words,
the condition $(*)$ is satisfied.
Thus in case that the condition $(*)$  fails to be satisfied,
there is a positive integer $e \geq 2$ such that,
for each $d$-element subset
$\{ x_{i_1}, x_{i_2}, \ldots, x_{i_d} \}$
of $\{ x_1, \ldots, x_n \}$
one has
$u = (x_{i_1} x_{i_2} \cdots x_{i_d})^e \in G(I)$.
Let $w = x_{n-d} x_{n-d+1}^{e-1} (\prod_{i=n-d+2}^{n} x_i^{e})
\in G(I)$.
Let $J$ denote the monomial ideal generated by those monomials
$v \in G(I)$ such that $v$ is bigger than $w$ with respect to the
reverse
lexicographic order.
Since $\prod_{i=n-d}^{n-1} x_i^{e} \in G(I)$,
by using Lemma \ref{juergen}
one has $w_0 = x_{n-d}w / x_n \in J$
and $w_1 = x_{n-d+1} w / x_n \in J$.
Thus the colon ideal $J : w$
is generated by a subset $M$ of $\{ x_{1}, \ldots, x_{n} \}$
with $\{ x_{1}, \ldots, x_{n-d}, x_{n-d+1}\} \subset M$.
Hence $q(I) \geq n - d + 1$, and 
thus we have $h(I) < q(I) + 1$, a contradiction.\hspace{9.2cm}
\end{proof}
 
As we pointed out in Section 1, a Cohen--Macaulay ideal is always unmixed. The converse is in general not true, even for matriodal ideals. For example, let $I \subset K[x_1, \cdots, x_6]$ be the monomial ideal
generated by $$x_1x_3, x_1x_4, x_1x_5, x_1x_6,
x_2x_3, x_2x_4, x_2x_5, x_2x_6,
x_3x_5, x_3x_6, x_4x_5, x_4x_6.$$
Then $I$ is matroidal and unmixed.  
However, $I$ is not Cohen--Macaulay.

It would, of course, be of great interest 
from a viewpoint of combinatorics
to classify all unmixed polymatroidal ideals.

\end{document}